\documentclass[12pt]{article}
\newtheorem{tm}{Theorem}
\newtheorem{Lemma}{Lemma}
\newcommand{\atsec}{\setcounter{equation}{0}\setcounter{Lemma}{0}}
\renewcommand{\thesection}{{\bf\S}\kern.3em\arabic{section}.\atsec}
\renewcommand{\thesubsection}{\kern.5em\arabic{section}.\arabic{subsection}.}
\renewcommand{\thesubsubsection}{\arabic{section}.\arabic{subsection}.\arabic{subsubsection}.}

\newcommand{\eps}{\varepsilon}

\newcommand{\eqbydef}{\stackrel{\rm def}{=}}

\newcommand{\dst}{\displaystyle}

\renewcommand{\theequation}{\arabic{section}.\arabic{equation}}

\renewcommand{\appendix}{
  \setcounter{section}{0}
  \renewcommand{\thesection}{\thesectionA}
  \renewcommand{\thesubsection}{\thesubsectionA}
  \renewcommand{\thesubsubsection}{\thesubsubsectionA}
  \renewcommand{\theequation}{\theequationA}
}
\newcommand{\thesectionA}{{\bf Appendix}\kern.3em\Alph{section}.\atsec}
\newcommand{\thesubsectionA}{\Alph{section}.\arabic{subsection}.}
\newcommand{\thesubsubsectionA}{\Alph{section}.\arabic{subsection}.\arabic{subsubsection}.}
\newcommand{\theequationA}{\Alph{section}.\arabic{equation}}
\def\Proof{\noindent{\sl  Proof}.$\,$}
\font\msam=msam10 scaled 1200
\font\Bbb=msbm10 scaled 1440

\newcommand{\Cc}{\mbox{\Bbb C}}
\newcommand{\Zz}{\mbox{\Bbb Z}}

\newcommand{\varkappa}{\mbox{\Bbb\symbol{"7B}}}

\newcommand{\BlackBox}{\hfill$\mbox{\msam \symbol{4}}$}
\newcommand{\WhiteBox}{\phantom{a}\hfill$\mbox{\msam\symbol{3}}$}

\newcommand{\Msk}{Moscow}
\def\rem#1{} 
\def\@listI{\leftmargin\leftmargini \parsep 1\p@
plus0.5\p@ minus0.2\p@ \topsep 2\p@ plus1\p@ minus1\p@
\itemsep 2\p@ plus0.5\p@ minus.5\p@}
\def\Binomial#1#2{\left({#1 \atop #2}\right)}
\def\Aa{B_{\alpha}}
\def\Alog{B'}
\def\Aas{\tilde B_{\alpha}}
\def\Alogs{\tilde B'}
\def\Ha{L_{\alpha}}
\def\Hlog{L'}
\def\Has{\tilde L_{\alpha}}
\def\Hlogs{\tilde L'}
\def\Ra{R_{\alpha}}
\def\Rlog{R'}
\def\DbL{K}
\def\Rs{D}
\def\Geg{P} 

\def\intone{\int\limits_{\kern-4pt -1}^1}
\def\intz{\int\limits_{0}^1}
\renewcommand{\Re}{\mathrm{Re}\,}
\renewcommand{\Im}{\mathrm{Im}\,}
\def\HGa#1#2#3#4{ \;{}_2 F_1\left(
  {{#1,\;\;#2} \atop {\;\;\;#3}}
  \,;\, #4
  \right) }
\def\HGu#1#2#3{ \;{}_2 F_1\left(
  {{#1,\;\;#2} \atop {\;\;\;#3}}
  \right) }
\def\HGb#1#2#3#4#5{ \;{}_3 F_2\left(
  {{#1,\;\;#2,\;\;#3} \atop {\;\;#4,\;\;#5}}
  \right) }
\def\HGc#1#2#3#4#5#6#7#8#9{ \;{}_5 F_4\left(
  {{#1,\;\,#2,\;\,#3,\;\,#4,\;\,#5} \atop
  {\;\,#6,\;\,#7,\;\,#8,\;\,#9}}
  \right) }
\def\sfrac#1#2{#1/#2}
\def\hfrac#1{(#1)/2}
\def\hodd#1{H^\mathrm{{odd}}_{#1}}

\def\hnum#1{H_{#1}}
\def\ubsubsection#1{{\bf #1}}
\begin{document}
\thispagestyle{empty}
\vbox{\large
\begin{center}
RUSSIAN ACADEMY OF SCIENCES\\[1ex]
M.V.KELDYSH INSTITUTE FOR APPLIED MATHEMATICS 
\end{center}
\vskip 2.5cm
\centerline{Preprint No 57/2003}
\vskip 3.5cm
\centerline{ S.Yu. Sadov}
\vskip 2cm
\centerline{COUPLING OF THE LEGENDRE POLYNOMIALS}
\vskip 2ex
\centerline{WITH KERNELS $\;|x-y|^{\alpha}\;$ and $\;\ln|x-y|\,$}
\vskip 7cm
\large
\centerline{\Msk, 2003}
}
\newpage\noindent
{\bf Abstract}

\medskip\noindent
Double integrals that represent matrix elements of the power and
logarithmic potentials (resp.\ $\,|x-y|^{\alpha}\,$ and $\,\ln|x-y|\,$)
in the Legendre polynomial basis are found in a closed form.
Several proofs are given, which involve different special functions
and identities. In particular,
a connection of the new formulae and Whipple's hypergeometric summation
formula is shown.

\bigskip
\noindent
{\bf E-mail:$\quad$ sadov@keldysh.ru}

\vspace{1cm}
\section{Introduction}
In this paper we evaluate the integrals
\begin{equation}
\label{powlegdef}
\Aa(m,n)=\intone\intone\,
|x-y|^{\alpha}\, P_m(x) \,P_n(y)\,dx\,dy
\qquad
(\Re\alpha>-1)
\end{equation}
(matrix elements of the Bessel potentials) $\,$ and
\begin{equation}
\label{loglegdef}
 \Alog(m,n)=\left.\frac{d}{d\alpha}\right|_{\alpha=0}\Aa(m,n)\,=
\intone\intone\,
\ln|x-y|\, P_m(x) \,P_n(y)\,dx\,dy\,.
\end{equation}
Here
$\,P_n
\,$ are the Legendre polynomials.
Notation and  relevant facts
regarding orthogonal polynomials and hypergeometric series are
collected in Appendix.

Sometimes it is more convenient to deal with the shifted
Legendre polynoimials (\ref{pdef}).
Define
\begin{equation}
\label{powlegdefs}
\Aas(m,n)=\intz\intz\,
|x-y|^{\alpha}\, p_m(x) \,p_n(y)\,dx\,dy
\qquad
(\Re\alpha>-1)
\end{equation}
and
\begin{equation}
\label{loglegdefs}
\Alogs(m,n)=\intz\intz\,
\ln|x-y|\, p_m(x) \,p_n(y)\,dx\,dy.
\end{equation}
It is readily seen that
$$
\begin{array}{l}
\Aa(m,n)\;=\;2^{\alpha+2}\;\Aas(m,n),\\[1ex]
\Alog(m,n)\;=\;4\,\left(\,\Alogs(m,n)\,+\,\delta_{m0}\,\delta_{n0}\,\ln 2\,
\right).
\end{array}
$$
Because of the symmetry of the integrands --- see
Appendix A.1.3. --- all the above integrals vanish when $\,m-n\,$ is odd.
In fact, we obtain more detailed results, which are non-trivial in the case
of odd $\,m-n\,$ as well, --- evaluation of the `halves' of the above
integrals. Instead of the interals over the square
$$
-1\,\leq x\,\leq\, 1,\quad
-1\,\leq y\,\leq\,1
\qquad(\,\mathrm{resp.}\;\;
0\,\leq x\,\leq\, 1,\quad
0\,\leq y\,\leq\,1 \,),
$$
consider the integrals of the same integrands over the triangle
$$
0\,\leq\,y\,\leq\,x\,\leq \,1 \qquad(\,\mathrm{resp.}\;\;
0\,\leq\,y\,\leq\,x\,\leq \,1\,),
$$
and denote them by
$\,\Ha(m,n)$, $\,\Hlog(m,n)$, $\,\Has(m,n)$, $\,\Hlogs(m,n)\,$
respectively.
Then
\begin{equation}
\label{Lsym}
\begin{array}{l}
 L(m,n)\,=\,(-1)^{m-n}\,L(n,m)\,;
 \\[1ex]
 B(m,n)\,=\,2\,L(m,n) \quad\mbox{\rm if}\;\;m-n\;\,\mbox{\rm is even},
\end{array}
\end{equation}
with $L$ being any of $\Ha$, $\Hlog$, $\Has$, $\Hlogs$ and the
corresponding $B$.
%
\noindent
Yet we formulate theorems on evaluation of the $B$'s separately,
because there are compact proofs that do not involve evaluation of
the $L$'s.

In view of the symmetry (\ref{Lsym}), it is enough to bring the results
assuming $\,m\geq n$. Denote $\,d=m-n$.  Set
\begin{equation}
\label{powlegrhs}
\Ra(m,n)\,=\,
 \frac{-1}{\alpha+1}\;
 \frac{(\hfrac{d-\alpha})_{n} \,(-\alpha-1)_{d}}%
  {(\hfrac{\alpha+d+4})_{n}\, (\alpha+2)_{d+1}}
\end{equation}
and
\begin{equation}
\label{loglegrhs}
\Rlog(m,n)\,=\,
 \dst\frac{1}{(m+n)(m+n+2)(d^2-1)}.
\end{equation}
Notice that if $\,d>0\,$ then the fraction $\Ra(m,n)$ can be reduced
by the common factor $\,(\alpha+1)$.

\begin{tm}
If $\,m-n\geq 0\,$ is even, then
\begin{equation}
\label{powlegs}
 \Aas(m,n)\;=\;2\,\Ra(m,n).
\end{equation}
\end{tm}

\begin{tm}
{\rm (i)}$\;$
If $\,m-n\,$ is even and $\,(m,n)\neq(0,0)$, then
\begin{equation}
\label{loglegs}
 \Alogs(m,n)\;=\;2\,\Rlog(m,n).
\end{equation}
{\rm(ii)}$\;$
Special case $\,m=n=0\,$:
\begin{equation}
\label{logleg0s}
 \Alogs(0,0)\, =\, -3/2.
\end{equation}
\end{tm}
Correspondingly, for even nonnegative $\,d=m-n$,
\begin{equation}
\label{powleg}
 \Aa(m,n)\,=\;2^{\alpha+3}\,\Ra(m,n),
\end{equation}
In particular, 
$$
 \Aa(0,0)=\frac{2^{\alpha+3}}{(\alpha+1)(\alpha+2)}.
$$
Also, assuming $\,m=n$ even and $\,(m,\,n)\neq(0,0)$, we have
\begin{equation}
\label{logleg}
 \Alog(m,n)\,=\;8\,\Rlog(m,n).
\end{equation}
The special case is
\begin{equation}
\label{logleg0}
 \Alog(0,0) = 4\ln 2-6 \;=\, -3.2274\dots\,.
\end{equation}

\begin{tm}
If $\,m-n\geq 0\,$ then
\begin{equation}
\label{powleghs}
 \Has(m,n)\;=\,\Ra(m,n).
\end{equation}
\end{tm}

\begin{tm}
{\rm(i)}$\;$
If $\,m-n\geq 2\,$ or $\,m=n\neq 0$, then
\begin{equation}
\label{logleghs}
 \Hlogs(m,n)\;=\,\Rlog(m,n).
\end{equation}
{\rm(ii)}$\;$ Special case $\,m-n=1\,$:
\begin{equation}
\label{logleg1s}
 \Hlogs(n+1,n)\;=\,-\,\frac{1}{q}\,\left(\hodd{2n+1}-\frac{1}{4}-
 \frac{1}{q}\right),
\end{equation}
where
$$
  q\,=\,(2n+1)(2n+3)\,=\,(m+n)(m+n+2),
$$
and $\,\hodd{2n+1}\,$ is the $\,n+1$-th odd harmonic sum,
\begin{equation}
\label{hodd}
 \hodd{2n+1}=\,1+\frac{1}{3}+\dots+\frac{1}{2n+1}.
\end{equation}
\end{tm}
The following numerical table supports Theorem 4 and includes the
omitted case $m=n=0$. Index $m$ enumerates the rows and $n$ the columns.
$$
 \Big[\Hlogs(m,n)\Big]_{m,n=0\dots 3}\;=\;\left[
 \begin{array}{cccc}
 \dst\frac{-3}{4}  &\dst\frac{5}{36} &\dst\frac{1}{24}   & \dst\frac{1}{120} \\[2.2ex]
 \dst\frac{-5}{36} &\dst\frac{-1}{8} &\dst\frac{61}{900} & \dst\frac{1}{72} \\[2.2ex]
 \dst\frac{1}{24}  &\dst\frac{-61}{900} &\dst\frac{-1}{24} & \dst\frac{527}{14700} \\[2.2ex]
 \dst\frac{-1}{120} &\dst\frac{1}{72} &\dst\frac{-527}{14700} & \dst\frac{1}{120}
 \end{array}
 \right]
$$

\medskip
The paper is organized as follows.

{\bf Section 2} contains a one-page proof of Theorem 2, due to
M.~Rahman \cite{Rah}.

{\bf Section 3}$\,$ contains a short proof of Theorem 1
proposed by R.~Askey. It involves three interesting formulae from
the theory of special functions. The proof is simple because nasty
details are worked out in the cited results.

In {\bf Section 4}$\,$ we prove Theorem 3 by
transforming the integral $\,\Has(m,n)\,$
to a certain terminating hypergeometric series,
which can be evaluated by a formula due to Whipple.
Moreover, through this connection we prove
the terminating case of Whipple's formula.
(Essense: our integral is a rational function with known zeros and
poles.)
A proof of Whipple's formula in the general case is then obtained
using a standard function-theoretic argument.

Finally, in {\bf Section 5}$\,$ we give an independent, self-contained,
and elementary proof of Theorem 4, which had been largely computer-aided.

\medskip
This work stems from my study of the integral $\,\Alog(m,n)\,$
that was originally motivated by numerical solution of
integral equations on polygons.
Formula (\ref{logleg}) was published without proof in
\cite{lviv99}, \cite{GSSS00}. 
The following problem about a generalization of the integral
$\,\Alogs(m,n)\,$ is important in those applications.

\medskip\noindent
{\bf Problem.}$\,$ {\it Propose an efficient numerical procedure
for evaluation of the integral
$$
\Alogs(m,n;\,\lambda)=
\intz\intz  \ln (\lambda x-y)\;
p_m(x)\, p_n(y) \, dx \,dy,
$$
where $\,\lambda\,$ is a complex parameter.
The method must be numerically stable and suitable for large
values of $m$ and $n$, and use the ordinary floating point arithmetics.}

\noindent
The $\,\Alogs(m,n;\,\lambda)\,$ are rational functions of $\lambda$,
if $m-n\geq 2$. Their evaluation using computer algebra with
exact rational arithmetics is not what we wish,
since as $|m-n|$ grows, the expressions become huge and the
calculation slow.

\bigskip
We conclude this section by deriving Theorem $2$ from Theorem $1$.

Case $\,n=d=0\,$ stands alone, but it is trivial.
Cases (i) $\,d\neq0\,$ and (ii) $\,d=0\,,n>0\,$ need separate
treatment, but in both cases $\,\Aas\,$ is a rational function of
$\,\alpha\,$ with a simple zero at $\,\alpha=0$. It remains to
factor out the $\,\alpha\,$ and to evaluate the remaining
quotient at $\alpha=0$. In case (i) the $\,\alpha\,$ is contained
in $\,(-\alpha-1)_d$. We have: as $\alpha\to 0$
$$
\frac{\Aas(m,n)}{\alpha}\,\to\,
\frac{2\,(\frac{d}{2})_n\,(1)_{d-2}}{(\frac{d}{2}+2)_n\,(2)_{d+1}}
\,=\frac{2\,\frac{d}{2}\,(\frac{d}{2}+1)}{%
(\frac{d}{2}+n+1)(\frac{d}{2}+n)(d-1)d(d+1)(d+2)}.
$$
Simplifying, we obtain (\ref{loglegs}).
In case (ii) the factor $\,\alpha\,$ comes from
$\,(d-\alpha/2)_n=(-\alpha/2)_n$. When $\alpha\to 0$,
$$
{\alpha}^{-1}\,{\Aas(n,n)}\;\to\;
{-(1)_{n-1}}/(2\cdot(2)_n)
\,=\,-1(2n(n+1))^{-1},
$$
in agreement with (\ref{loglegs}).
\section*{Acknowledgements}
Prof.\ {\it M.~Rahman} (Ottawa)
found the first proof of the formula (\ref{logleg})
and communicated it to me on Feb.~15, 1999.

Profs.\ {\it G.~Gasper} (Northwestern University), {\it R.~Askey}
(Wisconsin)  and
{\it E.D.~Krupnikov}
(Novosibirsk) responded to my question posted in an
online forum on special functions (opsftalk@nist.gov) in February 2001.
regarding formula (\ref{whip}); at that time I didn't know that
it was Whipple's.
Prof.\ Askey also provided a sketch of proof
of Theorem 1 (see Sect.~3).

\medskip
Publication of this work has been supported by the
Russian Foundation for Basic Research under grants
 02-01-01067 and 01-01-00517.

\section{M.~Rahman's proof of Theorem 2}
In this proof, the double integral is turned to
a sum of two ordinary integrals of a Legendre polynomial
and a Legendre function of the 2nd kind, for which the
values are known. We work with polynomials $P_n$, rather
than with the shifted Legendre polynomials, so the target
is formula (\ref{logleg}).

\begin{Lemma}
Denote
\begin{equation}
\label{rah2}
 F_n(x)=\frac{2n+1}{2}\,\int_{-1}^1 P_n(y)\,\ln |x-y|\, dy.
\end{equation}
Then
\begin{equation}
\label{rah3}
 F_n(x)\,=\,\frac{1}{2}\;{\rm P.V.}\intone
 \frac{P_{n+1}(y)-P_{n-1}(y)}{x-y}
 \,dy.
\end{equation}
\end{Lemma}
\Proof$\;$
Replacing $\,\ln|x-y|\,$ in (\ref{rah2}) by $\,\ln(x-y)$
doesnt't affect real part of the integral, whatever branch of the
logarithm is taken. Move the integration path to the
upper complex half-plane to obtain a regular
integrand, and integrate by parts, using (\ref{PRecd}).
The boundary terms vanish, since
$\;[P_{n+1}-P_{n-1}](\pm 1)=0$. $\,$
 Putting the integration path back on the real
line, we obtain (\ref{rah3}).  \WhiteBox

 \smallskip\noindent
The integral in (\ref{rah3}) is the difference of two Legendre
functions of the second kind
, $\,Q_{n-1}(x)\,$ and
$\,Q_n(x)\;$ \cite{Bat},$\,${\bf  3.6}$\,$(29).
Therefore,
\begin{equation}
\label{rah4}
\frac{2n+1}{2}\,\Alog(m,n) =
\int_{-1}^1 P_m(x)\, F_n(x)\,dx
= I(m,\,n-1) - I(m,\,n+1),
\end{equation}
where
$$
 I(m,k)\,=\,\int_{-1}^1  P_m(x)\, Q_k(x)\,dx.
$$
The integrals $\,I(m,k)\,$
are evaluated in \cite{Bat},$\,${\bf 3.12}$\,$(17). The result is:
\\[\smallskipamount]
(i)$\,$ if $\,m-k\,$ is even, then $\;\,I(m,k)=0\,$;
\\[\smallskipamount]
(ii)$\,$ if $\,m-k\,$ is odd, then
\begin{equation}
\label{Imk}
 I(m,k)\,=\,\frac{2}{(m-k)\,(k+m+1)}.
\end{equation}

\noindent
Assuming that $\,m-n=d\,$ is even,
substituion of (\ref{Imk}) to (\ref{rah4}) gives
$$
\frac{2n+1}{2}\,\Alog(m,n)\,=\,
\frac{2}{(d+1)\,(m+n)}\,
  -\,\frac{2}{(d-1)\,(m+n+2)},
$$
and (\ref{logleg}) easily follows.
\BlackBox

\section{R.~Askey's proof of Theorem 1}
In response to my question regarding formula (\ref{powleg})
posted in the online forum on special functions, R.~Askey
\cite{Asopsf} wrote:
\footnote{Citation is slightly abridged.}

\begin{sl}{
 There is a long paper by Polya and Szego in the early 1930s,
see either Polya's collected papers or those of Szeg\"o,
in which
they find an expansion of $\,|x-y|^{\alpha}\,$ as a sum of products of
ultraspherical polynomials.
Then the integral (\ref{powlegdef}) can be evaluated as a single
hypergeometric series by using the formula for the integral of an
ultraspherical polynomial times a Legendre polynomial.
This is a special case of a formula of Gegenbauer which I have used
frequently.  There are a number of relatively simple derivations of this
formula.  One of the easiest is the one Gasper and I gave in the paper
we wrote for the meeting celebrating de Branges's proof of the
Bieberbach conjecture.  It just uses the generating function of
ultraspherical polynomials, differentiation, and orthogonality to get
the result.
  Then it is just a matter of seeing if the hypergeometric
series which comes from these calculations can be summed.
Sadov claimed a specific formula is true, so the series can be summed.
}\end{sl}

This section supplies details of Askey's approach.

\noindent
{\bf Step 1. Decomposition.}$\;$
Ultraspherical (or Gegenbauer's) polynomials
$\;\Geg^{(\nu)}_n(x)\;$ are defined by the generating function
\begin{equation}
\label{gfgeg}
(1\,-\,2\,x\,t\,+t^2)^{-\nu}\;=\;
\sum_{n=0}^{\infty} \Geg^{(\nu)}_n(x) \;t^n.
\end{equation}
The decomposition of $\,|x-y|^{\alpha}\,$ due to Polya and Szeg\"o
is $\,$
\cite[S.~27]{PS31} 
\begin{equation}
\label{PolSze}
|x-y|^{\alpha}
\;=\;M_1
\,\sum_{j=0}^{\infty}\left(1+\,\frac{j}{\nu}\right)
\;\Geg_j^{(\nu)}(x) \;\Geg_j^{(\nu)}(y),
\end{equation}
where $\,\nu=-\alpha/2\;$  and
\begin{equation}
\label{PSprefac1}
M_1\,=\,M_1(\alpha)\,=\,\frac{\Gamma\left(\frac{1}{2}-\nu\right)\,
\Gamma(1+\nu)}{\Gamma(\frac{1}{2})}.
\end{equation}
In (\ref{PolSze}), $\,\alpha>0$. 
Convergence in the r.h.s.\ follows from the estimation
$$
 \Geg_j^{(\nu)}(x)=O(j^{\nu-1}), \quad\;-1\leq x\leq 1.
$$

\noindent
{\bf Step 2. Integration.}$\;$
The formula due to Gegenbauer reads \cite{AsGas}
\begin{equation}
\label{gegleg}
\intone \Geg_j^{(\nu)}(x)\,P_{j-2s}(x)\,dx\;=
\;\frac{(\nu)_{j-s}\,(\nu-1/2)_s}
{s!\,(1/2)_{j-s+1}}.
\end{equation}
Assume $\,m-n=d\geq 0\,$ and even.
Integral
of the $j$-th term in (\ref{PolSze}) times
$\,P_m(x)\,P_n(y)\,$ may be nonzero only when $\,j=m+2k$,
$\,k\geq 0$.
We have
%
$$
\begin{array}{ll}
\dst
\langle
\Geg_j^{(\nu)}(x),\,P_m(x) \rangle &=\dst
\;\frac{(\nu)_{m+k}\,(\nu-1/2)_k}
{k!\,(1/2)_{m+k+1}},
\\[3ex]
\langle
\Geg_j^{(\nu)}(y),\,P_n(y) \rangle &=\dst
\;\frac{(\nu)_{n+k+d/2}\,(\nu-1/2)_{k+d/2}}
{(k+d/2)!\,(1/2)_{n+k+d/2+1}}.
\end{array}
$$
Write the obtained sum in the standard hypergeometric form
(\ref{hgdef}), with a prefactor, using the identity
$$
1- \frac{2\,j}{\alpha}\,=\left(1+\frac{m}{\nu}\right)\;
\frac{((m+\nu)/2+1)_k}{((m+\nu)/2)_k},
$$
and the identities of the type (\ref{Pohsum}), for example,
$$
(\nu)_{n+k+d/2}\,=\,(\nu)_{n+d/2}\,(\nu+n+d/2)_{k}.
$$
Also recall that $\,n+d/2=m-d/2$.
The result of integration takes the form
\begin{equation}
\label{PSint}
\Aa(m,n)\,=
\,M_1\,
M_2\,\cdot\,
\HGc{\nu+m}{1+\frac{m+\nu}{2}}{\nu-\frac12}{\nu+m-\frac{d}{2}}%
{\nu+\frac{d-1}{2}}{\frac{m+\nu}{2}}{\frac32+m}{\frac{d}{2}+1}%
{\frac{3-d}{2}+m}
\end{equation}
with
\begin{equation}
\label{PSprefac2}
M_2\,=
\,\frac{(1+m/\nu)\,
(\nu)_m\,(\nu)_{m-d/2}\,(\nu-\sfrac{1}{2})_{d/2}}%
{(\sfrac{1}{2})_{m+1}\;(1)_{d/2}\;(\sfrac{1}{2})_{m-d/2+1}}.
\end{equation}

\noindent
{\bf Step 3. Summation.}$\;$
The above series
$_5F_4$ can be summed by the limiting case
(\ref{dougall54})
of Dougall's formula for a well-poised terminating $_7F_6$.
The result is a ratio of products of $\Gamma$ functions.
Arguments of four $\,\Gamma$'s in the numerator are
\begin{equation}
\label{dougnum}
m+3/2, \quad
d/2+1,\quad
m+(3-d)/2, \quad
\alpha+2.
\end{equation}
In the denominator we have four $\Gamma$'s with arguments
\begin{equation}
\label{dougden}
m+1-\alpha/2, \quad
(3+\alpha)/2,\quad
m+2+(\alpha-d)/2, \quad
(3+\alpha+d)/2.
\end{equation}
In total, we obtain an expression for $\Aa(m,n)$ in the form
of a fraction with 10 factors in the numerator and 8 factors in the
denominator. To simplify, we break down the triple product
(\ref{PSint}) into 5 groups and write
$$
\Aa(m,n)\;=\;\Pi_1\;\Pi_2\;\Pi_3\;\Pi_4\;\Pi_5.
$$
In the expressions for the $\Pi$'s below, we show
which little factors come from which big ones.
The notation [$j$,\,$k$] stands for the
$\,k$-th term in the numerator or denominator
of the $j$-th big factor: $\,j=1\,$ refers to
$\,M_1$; $\,j=2$ to $M_2$; and $\,j=3$ to $_5F_4$.
Denominator terms are shown in bold.
For example, [2,1] is the first term in the numerator of
(\ref{PSprefac2}), that is, $\,(1+m/\nu)$.
One cancelation $\;{\rm [3,2]}/\mathbf{[2,2]}=1\;$ occurs immediately
and those terms are not attributed to any group.
The grouping of other terms and the results of simplification are
as follows
$$
\begin{array}{l}
\dst
\Pi_1\,=\,\sfrac{\mathrm{[1,2]\;[2,1]\;[2,2]}}{\mathbf{[3,1]}}\;=\;1,
\\[1ex]
\dst
\Pi_2\,=\,\sfrac{\mathrm{[3,1]}}{\mathbf{[1,1]\;[2,1]}}\;=\;
\Gamma(m+3/2),
\\[1ex]
\dst
\Pi_3\,=\,\sfrac{\mathrm{[1,1]}}{\mathbf{[3,2]}}\;=\;
2\,/(\alpha+1),
\\[1.6ex]
\dst
\Pi_4\,=\,\sfrac{\mathrm{[3,4]}}{\mathbf{[3,3]\;[3,4]}}\;=\;
2^{2+\alpha+d}\,\pi^{-1/2}\;\left[\,(\alpha+2)_{d+1}\;(2+
{\textstyle\frac{\alpha+d}{2}})_{n}\right]^{-1}
\\[2.1ex]
\dst
\Pi_5\,=\,\frac{\mathrm{[2,2]\;[2,3]\;[3,3]}}{\mathbf{[2,3]}}\;=\;
\frac{(-\alpha/2)_{m-d/2}\,(-(\alpha+1)/2)_{d/2}\,\Gamma(m+
\frac{3-d}{2})}{(1/2)_{m-d/2+1}}.
\end{array}
$$
The part $\;\Gamma(\dots)/(1/2)_{m-d/2+1}\;$ in $\,\Pi_5\,$
cancels with $\,\pi^{-1/2}\,$ from $\,\Pi_4$. It remains
to apply the duplication formula (\ref{Pohdup}) with $\,a=-\alpha/2\,$,
$\;n=d/2$, and (\ref{powleg}) follows.
\BlackBox

\subsection*{Limiting case: another proof of Theorem 2}
A proof of Theorem 2, different from those given in Sections 2 and 5,
can be obtained as a special limiting case of the above argument.
It follows the same lines, but instead of the ultraspherical
polynomials we encounter Chebyshev's polynomials $T_n$,
where analogues of (\ref{PolSze}) and (\ref{gegleg}) are rather elementary.
Chebyshev's polynomials show up here, because
$\;\left.\partial_{\nu}\,\Geg_n^{(\nu)}(x)\right|_{\nu=0}
\;=({2}/{n})\,T_n(x)$ $\; (n>0)$
and the generating function for $T_n$'s is
$$
\ln\,(1\,-\,2\,x\,t\,+t^2)\;=\;-2\,
\sum_{n=1}^{\infty} n^{-1}\, T_n(x) \;t^n.
$$
The role of (\ref{PolSze}) is played by the formula
$$
 \ln |x-y|\;=\;-\ln 2\,-\,2\,
 \sum_{j=1}^{\infty} j^{-1}\,T_j(x)\,T_j(y),
$$
which can be identified with the well known Fourier expansion
$$
-\ln|2\,\sin (\theta/2)|\,=\,\sum_{j=1}^{\infty} |j|^{-1}\,\cos(j\theta)
$$
using the definition of Chebyshev's polynomials
$\;T_j(\cos\theta)=\cos(j\theta)$.

\noindent
The analogue of (\ref{gegleg}) reads $\,$
(for $\,j>0$)
$$
\frac{2}{j}\;\,\intone T_{j}(x)\,P_{j-2s}(x)\,dx\;=\;-\,
\frac{(s+1)_{j-2s-1}}{2\,(s-1/2)_{j-2s+2}}.
$$
Step 3 again amounts to evaluation of the particular case of
the series (\ref{PSprefac2}) with $\,\nu=0\,$
using Dougall's formula. (It seems attractive to find
its elementary evaluation in this case.)
\section{Proof of Theorem 3 via Whipple's formula}

\subsection{Reduction of integral $\,\Has(m,n)\,$ to a terminating
$_3F_2$ series}
Set
\begin{equation}
\label{ialpha}
A_{\alpha}(n,\,x) = \int
_0^x (x-y)^{\alpha} p_n(y)\,dy,
\end{equation}
then
\begin{equation}
\label{Aat}
\Has(m,n)\,=\,
\int
_0^1 p_m(x) \,A_{\alpha}(n,\,x)\,dx.
\end{equation}
\begin{Lemma}
\begin{equation}
\label{basx}
 A_{\alpha}(n,\,x)\,
 =\,
 \frac{x^{\alpha+1}}{\alpha+1} \,
 \HGa{-n}{n+1}{\alpha+2}{x}.
\end{equation}
In particular, 
\begin{equation}
\label{bas1}
 A_{\alpha}(n,\,1)  \,
 =\, \frac{(-1)^n\, (-\alpha)_n}{(\alpha+1)_{n+1}}.
\end{equation}
\end{Lemma}
{\Proof}$\;$
Apply the Rodrigues formula (\ref{pRod}) and integrate
by parts $n$ times.
Assume $\,\alpha>n$, then the
boundary terms vanish and we get
$$
 A_{\alpha}(n,\,x)
 \,=\,\int\limits%
 _0^x \frac{(-\alpha)_n \,(-1)^n}{n!}\, (x-y)^{\alpha-n}
 \;y^n\,(1-y)^n\,dy.
$$
Setting $\,y=tx$ leads to the integral (\ref{f21int}) for
$_2F_1$ with parameters
$\, a=-n$, $\,b=n+1$, $\,c=\alpha+2$,
and the prefactor
$$ x^{\alpha+1}\;
 \frac{(-1)^n\,
(-\alpha)_n}{n!}\;\frac{\Gamma(n+1)\Gamma(\alpha-n+1)}{\Gamma(\alpha+2)},
$$
which reduces to $\,x^{\alpha+1}/(\alpha+1)\,$. 
The assumption $\,\alpha>n\,$ can be dropped because both parts
of (\ref{basx}) are analytic in $\,\alpha$.
\WhiteBox

\begin{Lemma}
\begin{equation}
\label{ahg}
\Has(m,n)\;=\;\frac{(-\alpha-1)_{m}}{(\alpha+1)_{m+2}}
\;\HGb{-n}{n+1}{\alpha+2}{\alpha+m+3}{\alpha+2-m}.
\end{equation}
\end{Lemma}
\Proof$\;$
Substituting (\ref{basx}) in (\ref{Aat}), we find
$$
\Has(m,n)\,=\,\sum_{j=0}^{n}\frac{(-n)_j\, (n+1)_j}{j!\,(\alpha+1)_{j+1}}\;
\int\limits%
_0^1\, x^{\alpha+1+j}\,p_m(x)\,dx.
$$
Since $\,p_m(x)=(-1)^m\,p_m(1-x)$,
the integral in the common term is
$$
 (-1)^m\,
 I(m,\,\alpha+1+j,\,1)\,=\,
 \frac{(-\alpha-1-j)_m}{(\alpha+2+j)_{m+1}},
$$
according to (\ref{bas1}).
Using the identities (cf. (\ref{Pohsum}))
$$
 (\alpha+1)_{j+1}\, (\alpha+j+2)_{m+1}\,=\,
 (\alpha+1)_{m+2}\,(\alpha+m+3)_{j}
$$
and
$$
(-\alpha-j-1)_m\,=\,(-\alpha-1)_m \frac{(\alpha+2)_j}{(\alpha-m+2)_j},
$$
the common term is transformed into
$$
 \frac{(-\alpha-1)_m\,(-n)_j\,(n+1)_j\,(\alpha+2)_j}%
 {(\alpha+1)_{m+2}\,(\alpha+m+3)_{j}\,(\alpha+2-m)_{j}},
$$
so (\ref{ahg}) follows.
\WhiteBox

\medskip\noindent
{\bf Remark.}
The symmetry (\ref{Lsym}) apparent in the l.h.s.\
of (\ref{ahg}) is hidden in the r.h.s., where it
is a particular case of Thomae's formula (\ref{thomae32}).

\subsection{\bf Equivalence of Theorem 1 and the terminating
case \\
of Whipple's formula}
Invoke the terminating case of Whipple's summation formula
(\ref{whipterm}) with
\begin{equation}
\label{subwhip}
\;b=\alpha+m+3, \quad\;\;c=\alpha+2-m,
\end{equation}
to evaluate the r.h.s.\ of (\ref{ahg}), and get
$$
 \Has(m,n)=4^{n}\,\frac{(-\alpha-1)_m\,(\frac{\alpha+d+3}{2})_n\,
 (\frac{\alpha+2-m-n}{2})_n}
 {(\alpha+1)_{m+2}\,(\alpha+m+3)_n\,(\alpha+2-m)_n}.
$$
To identify this expression with (\ref{powlegrhs}),
decompose the ratio of the right sides of the two formulae
into three groups and simplify using (\ref{Pohsum})--(\ref{Pohdup}):
$$
\begin{array}{c}
\dst
\frac{4^{n}\,((\alpha+d+3)/{2})_n\,((\alpha+d+4)/{2})_n}%
{(\alpha+1)_{m+2}\,(\alpha+m+3)_n\,(\alpha+d+3)_n}
\;=\,\frac{1}{\alpha+1},
\\[3ex]\dst
\frac{(-\alpha-1)_m}{(-\alpha-1)_d\,(\alpha+2-m)_n}\;=\,
\frac{((\alpha+2-m-n)/{2})_n}{((d-\alpha)/{2})_{n}}
\;=\,(-1)^n.
\end{array}
$$
Conversely, (\ref{powleg}) implies (\ref{whipterm})
in the case when the parameters $\,b\,$ and $\,c\,$
can be presented in form (\ref{subwhip}) with $\,m\in\Zz_+\,$
and $\,\Re\alpha>-1$, that is if $\,b-c\,$ is an even non-negative integer
and $\,\Re(b+c)>4$. But both parts of (\ref{whipterm}) are rational
functions of $\,b\,$ and $\,c\,$. Hence the restrictions
can be dropped.
\BlackBox

\subsection{Proof of Whipple's formula in the terminating case}
The l.h.s.\ of (\ref{whipterm}) is a rational function of $\,\alpha\,$
of the form
\begin{equation}
\label{Q}
\HGb{-n}{n+1}{\alpha+2}{\alpha+m+3}{\alpha+2-m}=
\frac{Q(\alpha)}{(\alpha+m+3)_n (\alpha-m+2)_n},
\end{equation}
where $\,Q(\alpha)\,$ is a polynomial of degree $\,2n\,$
(of course $Q$ depends on $m$ and $n$).
We shall find all roots of $\,Q(\alpha)\,$ by
using relation (\ref{ahg}) and studying zeros and
poles of $\,\Has(m,n)$. In this way we obtain a self-contained
proof of Theorem 1 together with formula (\ref{whipterm}).
%
\begin{Lemma}
If $\,m-n\,>0\,$ and even. Then $\,$
{\rm (i)} $\,$
$\,\Has(m,n)=0\;$ for $\,\alpha=0,\,2,\dots, m+n-2\,;$
$\;$ {\rm (ii)} $\,$
$\,Q(m+n-2s)=0\,$ for $\,s=1,\dots,n$.
\end{Lemma}

\Proof$\;$ (i)$\,$
It is enough to show that $\,\Aas(m,n)=0\,$ at the above
values of $\alpha$. In this case,
$
 \,|x-y|^{\alpha}\,
$
is a polynomial of degree $\,<m+n$.
In its binomial expansion $\,c_{kl}\,x^k\,y^l$,
the exponents $k$, $l$ satisfy at least one of the
inequalities $\,k<m$, $\,l<n\,$. Therefore, all products
$\,\langle x^k,\,p_{m}(x)\rangle\,\langle y^l,\,p_{n}(y)\rangle\,$
are $0$.

\noindent
(ii)$\,$
Using the relation
$$
\sfrac{(-1-\alpha)_m}{(\alpha+2-m)_n}\,
=\,(-1)^n\,(-\alpha-1)_{m-n},
$$
we get
$$
 \Has(m,n)\,=\, R(\alpha)\, Q(\alpha) \;
 \alpha\,(\alpha+1)\dots(\alpha+m-n-2),
$$
where $R(\alpha)$ is a rational function, which does not have poles
and zeros at
nonnegative integers.
Now the claim follows from (i).
\WhiteBox

In the next Lemma, we find $\,n\,$ other roots of $\,Q(\alpha)$.

\begin{Lemma}
Let $\,m>n$. Then $\,${\rm (i)}$\,$
$\,\Has(m,n)\,$ does not have poles at
$\,\alpha=-m-n-1+2s$, $\,s=0,\dots,n-1\,;$
$\;${\rm (ii)}$\,$
$\;Q(-m-n-1+2s)=0\;$ for
$\,s=0,\dots,n-1$.
\end{Lemma}
\Proof$\;$
(ii) follows from (i), because the
values $\,\alpha=-m-n-1+2s$, $\;0\leq s\leq n-1$,
are uncompensated
roots of the denominator in the r.h.s.\ in (\ref{ahg}),
while $\Has(m,n)$ is regular at these values.

\noindent
To prove (i),  consider analytical continuation
of $\Has$. Using the identity
$$
r^\alpha=(\Gamma(-\alpha))^{-1}\,\int
_{0}^{\infty} e^{-rt}\, t^{-\alpha-1}\,dt
\qquad
(\Re\,\alpha<0),
$$
we obtain for $\;-1<\Re\,\alpha<0$
$$
\Has(m,n)=\frac{1}{\Gamma(-\alpha)}\;\int\limits_0^{\infty}
t^{-\alpha-1}\;\int\limits_{0}^1 e^{-tx} p_m(x)\,
\int\limits_0^x e^{ty} p_n(y) \,dy\,dx\,dt.
$$
Divide the $t$-domain into
$\int_0^1$ and $\int_1^{\infty}$; the first part is a regular
function of $\,\alpha\,$ for $\,\Re\alpha<0$.
Let us find asymptotics of the integrand as $\,t\to\infty\,$
in the second part and study its analytical continuation
in $\,\Re\,\alpha\leq -1$.

The following formula of indefinite
integration
\begin{equation}
\label{expp}
\int p(\xi) e^{\lambda\xi}\,d\xi=
e^{\lambda\xi} (\lambda^{-1} p(\xi) -
\lambda^{-2} p'(\xi) + \lambda^{-3} p''(\xi) -\dots)\,
\end{equation}
holds for any polynomial $p(\xi)$. So we have
$$
\int\limits_0^x e^{ty} p_n(y)\,dy=
e^{tx}\sum_{j=0}^n \frac{(-1)^j}{t^{j+1}}\, p_n^{(j)}(x)
\,-\,\sum_{j=0}^n \frac{(-1)^j}{t^{j+1}}\, p_n^{(j)}(0).
$$
Applying (\ref{expp}) once again and using orthogonality,
we get
$$
\int\limits_0^1 e^{-tx}\, p_m(x) \,\int\limits_0^x e^{ty} \,p_n(y)\,dy\,dx
\;=\;
\sum_{j=0}^n\sum_{k=0}^m \frac{(-1)^{j+1}}{t^{j+k+2}}\; p_n^{(j)}(0)\,
p_m^{(k)}(0) \,+\, O(e^{-t}).
$$
The coefficient at $\,t^{-r}\,$
is the residual of $\,\Has(m,n)\,$ at the pole
$\alpha=-r$,  due to
$$
 \int_1^{\infty} t^{-\alpha-1-r}\,dt=\frac{1}{\alpha+r}.
$$
By (\ref{pMacid}), these coefficients vanish
for $\,r=m+n+1-2s$, $0\leq s\leq n-1$. Hence there are no poles
at the values stated in (i) of the Lemma.
\WhiteBox

\bigskip\noindent
The two Lemmas supply the complete list
of $2n$ roots of the numerator $Q(x)$. Thus we obtain:
$\;$
{\it If $\,m-n$ is even
and $>0$, then
\begin{equation}
\label{f32c}
\HGb{-n}{n+1}{\alpha+2}{\alpha+m+3}{\alpha+2-m}\;=\;
C\; \frac{(\frac{\alpha+m+3-n}{2})_n\, (\frac{\alpha-m+2-n}{2})_n }
{(\alpha+m+3)_n\, (\alpha-m+2)_n}\, ,
\end{equation}
where $C$ is independent of $\alpha$.}

The left side is also a rational function of
$\,m\,$ with degrees $2n$ of the numerator and denominator,
so the coefficient $C$ is independent of $m$.
Moreover, we can dismiss the above assumption of $m-n$
being even and positive and allow $m$ to be non-integral.
To determine $C$ (which now may depend only on $n$), let us take
$\,\alpha\to\infty$. Then the l.h.s.\ of
(\ref{f32c}) is $\;1+o(1)\,$, while the r.h.s.\ is $\;C/4^n\,+\,o(1)$.
Thus $\;C=4^n\,$ and (\ref{whipterm}) follows by substitution
(\ref{subwhip}).
\BlackBox

\subsection{Proof of Whipple's formula in the general case}
Fix the parameters $b$ and $c$ in Whipple's formula (\ref{whip}).
The difference of the l.h.s.\ and the r.h.s.\
of (\ref{whip}) is an entire function $\,f(a)$.
The motif of the proof is simple
(cf.\ technique based on Carlson's theorem \cite[Chapter V]{Bai}).
The function $\,f(a)\,$ has sufficiently many known zeros
and a uniformly controlled growth as
$\,|a|\to \infty$.
From this we will conclude that $\,f(a)\equiv 0$.

At this time, the proof is not self-contained, though marginally.
I have to recourse to additional facts from the theory of generalized
hypergeometric series: Karlsson's formula and Dixon's formula.

\begin{Lemma}$\,$
{\rm (i)}$\,$ If $\,\Re(b+c)>1\,$ then
the right-hand side of (\ref{whip}) has estimation
$\;O(e^{\pi |\Im a|} \,(1+|a|)^{1-\Re(b+c)})\,$  as $\,|a|\to\infty$.
\\[\smallskipamount]
{\rm (ii)}$\,$ If $\,b$, $c\,$ are real and  $\,b+c>1\,$ then
the left-hand side of (\ref{whip}) is
$\;O(e^{\pi|a|} \,(1+|a|)^{1-(b+c)})\,$
as $\,|a|\to\infty$.
\end{Lemma}

\Proof$\;$
(i) $\,$ follows from the inequality
\begin{equation}
\label{invgam}
({\Gamma(z)\Gamma(r-z)})^{-1}
\;<\; C(r)\,e^{\pi|\Im z|}\,(1+|z|)^{1-{\Re r}},
\end{equation}
which is a consequence of the formula $\Gamma(z)\Gamma(1-z)=
\pi/\sin(\pi z)$
and the estimation
$\;|\Gamma(x)/\Gamma(x+r)|=O(x^{-\Re{r}}),$
$\,x\,\to\infty$, $|\arg x|<\pi$.

\smallskip\noindent
(ii) $\,$ The series has the form
$$
 \sum_{n=0}^{\infty} \mu_n\;\frac{(a)_n\,(1-a)_n}{n!\,((b+c+1)/2)_n},
$$
where
$$
\mu_n=\frac{((b+c-1)/2)_n\,((b+c+1)/2)_n}{(b)_n\,(c)_n}\,=O(1)
\quad\;\mbox{\rm as}\;\;n\to\infty.
$$
Next,
$\,(a)_n\,(1-a)_n\,$
is a series with positive coefficients in $\,t=(a-1/2)^2$.
So maximum of our $_3F_2$ on the circle $\,|a-1/2|=R\,$ is bounded,
up to a constant, by the respective maximum of
$$
 \HGu{a}{1-a}{(b+c+1)/2}
 \;\;\stackrel{\mbox{\rm(\ref{hgGauss})}}{=}\;\;
 \frac{\Gamma(\frac{b+c+1}{2})\,
 \Gamma(\frac{b+c-1}{2})}{\Gamma(\frac{b+c+1}{2}-a)\,
 \Gamma(\frac{b+c-1}{2}+a)}.
$$
Finally, by (\ref{invgam}), the r.h.s.\ is
$\,O(e^{\pi|\Im a|}\,(1+|a|)^{1-(b+c)})$.
\WhiteBox

\smallskip 
If in (ii) we had $\,e^{\pi|\Im a|}\,$ instead of
$\,e^{\pi|a|}\,$, then the proof could be easily completed
via the maximum principle applied to
the function $\,f(a)\,\mathrm{cosec}\,\pi a$, which is entire.
(No poles on $\Zz$: the terminating case is already established).

In reality, we have only found that
 $\,\ln |f(a)| \leq \,(\pi+o(1)) |a|
 $.
Consider the asymptotical density of zeros of $f$,
$$
\nu=\,\overline{\lim\limits_{R\to\infty}}
\,N(R)/R,
\qquad
N(R)=\#\{\,a\;:\,\; f(a)=0\;\,\mbox{\rm and}
\;\,|a|<R\,\}
$$
It follows from Jensen's formula that the estimation $\,\nu>2\,$
would imply $\,f\equiv 0$.
Thus far we only know that $f|\Zz=0$,
which gives $\,\nu\geq 2\,$ but leaves a possibility of the equality.
To complete the proof, we shall present another infinite sequence of zeros.
Notice that the r.h.s.\ of (\ref{whip}) vanishes
when $\,a=-b-2k$, $\,k\in\Zz$.
Show that those are also zeros of the l.h.s. (By symmetry between
$b$ and $c$, yet more zeros will be found.)

\begin{Lemma}
Let $n$ be a nonnegative even integer, and $\,b$, $c\,$ arbitrary, except
nonnegative integers, and denote $\,d=(b+c-1)/2$.  Then
\begin{equation}
\label{spcase2}
 \HGb{-b-n}{b+n+1}{d}{b}{c}=0.
\end{equation}
\end{Lemma}
\Proof$\;$
Through Karlsson's formula (\ref{Karlsson}) and
Gauss' formula (\ref{hgGauss}) we transform (\ref{spcase2}) into
the terminating $_3F_2$ with parameters
$\,({-n-1},\,{-b-n},\,{d};$ ${b},\,{-d-n})$,
times a nonsingular factor.
Such an $\,{}_3 F_2\,$ falls under
Dixon's formula 
\cite[3.1]{Bai},
\cite[4.4 (5)]{Bat}, the sum being
a ratio of products of $\Gamma$ functions.
The only singular factor,
$\Gamma(-n)$, appears in the denominator.
(The factor $\,\Gamma((1-n)/2)\,$ is safe,
because $n$ is even.)
Therefore, the ratio equals to $0$.
\BlackBox


\section{Proof of Theorem 4}
\label{Sectlog2}
The double integral $\,\Hlogs(m,n)\,$
will be expressed via a double sum
using the Maclaurin expansions of the shifted Legendre polynomials.
We rearrange the double sum and show that new partial single sums
satisfy certain recurrence relation, from which the theorem follows.

We set aside the simple case $\,m=n=0\,$ and
prove formulae (\ref{logleghs}) and (\ref{logleg1s}) assuming
$\;m=\max(m,n)>0$.

\subsection{Reduction of the integral $\,\Alogs\,$ to a finite
sum}

Below {\it any function} means a function as good as a
polynomial.

Denote the average value of a function $f$ on an interval $\,(a,b)\,$ by
$$
\langle f\rangle_{a}^{b} = (b-a)^{-1}\,\int_a^b f(t)\,dt.
$$
Introduce two operators, $M$ and $N$, as follows:
\begin{equation}
\label{opM}
  Mf(x)= \langle f\rangle_{0}^{x},
\end{equation}
\begin{equation}
\label{opN}
  Nf(x)=\int_0^x \langle f\rangle_{t}^{x}\,dt.
\end{equation}
Clearly, if $f$ is a polynomial of degree $n$ then $\,Mf\,$ is a polynomial
of degree $\,n\,$, and $\,Nf\,$ a polynomial of degree $\,n+1$.

\begin{Lemma}
Let $f$, $g$ be any functions, and
$\;G(x)=\dst\int_0^x g(y) dy$. Then
\begin{equation}
\label{getridlog}
 \int\limits_0^1\,f(x)\,\int\limits_0^x \ln(x-y)\, g(y)\,dy\,dx=
 -\int\limits_0^1 M[f\cdot G](x)\,dx -\int\limits_0^1 f(x)\cdot Ng(x)\,dx.
\end{equation}
\end{Lemma}

\Proof$\;$ Integration by parts w.r.t.\ $y\,$ yields
$$
\int_0^x \ln(x-y)\, g(y)\,dy\;=\;
\ln x\;G(x)\,-\,Ng(x).
$$
The assertion of Lemma follows by second integration
$$
 \int_0^1 \ln x\;f(x)\, G(x)\,dx\, =-\int_0^1 M[f\cdot G](x)\,dx,
$$
where we remember that $\,\lim_{x\to 0} G(x)\ln x=0$.
\WhiteBox

\medskip
\noindent
{\bf Definition of the sums $\,S(m,n,r)$}.$\;$
Let $0\leq r\leq (m+n)$.
The quantities $\,S(m,n,r)\,$ are the coefficients of the generating
function
\begin{equation}
\label{gfSmnp}
p_m(x)\, q_{n+1}(x)\;=\,\sum_{r=0}^{m+n} S(m,n,r)\, x^{r+1},
\end{equation}
where $\,q_{n+1}\,$ are the integrated shifted Legendre polynomials
--- see (\ref{ipMac}).
The explicit formula is $\,$ (cf.\
(\ref{pMac}) and (\ref{ipMac}) $\,$)
\begin{equation}
\label{rec1}
  S(m,n,r)\,=\,(-1)^{m+n-r}\,\sum_{j=j{\rm min}}^{j{\rm max}}\,
  \frac{(m+j)!}{(m-j)!\,j!^2}\;\frac{(n+k)!}{(n-k)!\,k!\,(k+1)!}\,,
\end{equation}
where
$\;k=r-j$, $ \;j{\rm min}=\max(0,r-n)$, $\;j{\rm max}=\min(m,r)$.

\noindent
The promised double sum
and its rearrangement  are as follows
\begin{equation}
\label{dblsum}
\begin{array}{c}
\dst
\sum_{j=0}^m \sum_{k=0}^n \frac{(-1)^{m-j}\,(m+j)!}{(m-j)!\,j!^2}\;
\frac{(-1)^{n-k}\,(n+k)!}{(n-k)!\,k!\,(k+1)!}\;\frac{1}{(j+k+2)^2}
\\[3ex]
=\;
\dst\sum_{r=0}^{m+n}\frac{S(m,n,r)}{(r+2)^2}
\;\eqbydef\;\DbL(m,n).
\end{array}
\end{equation}
\begin{Lemma}
\label{logred}
For any $\,m$, $n$
\begin{equation}
\label{ilog2sum}
 \int\kern-1ex\int\limits_{\kern-2ex 0<y<x<1}
  p_m(x)\,p_n(y)\,\ln(x-y)\,dy\,dx\,=\,
 -\DbL(m,n)\,-\,
  \langle p_m,\,Np_n\rangle.
\end{equation}
\end{Lemma}

\Proof$\;$
Evaluate the functional
$\;u\to\langle Mu,\,1\rangle\,$,
which appears in (\ref{getridlog}),
in terms of the Maclaurin expansion of the function $\,u(x)$
\begin{equation}
\label{genericMac}
 u(x)=\sum\nolimits_{r\geq 0}\, u_r\,x^r.
\end{equation}
The operator $M$ divides the $r$-th coefficient by $\,r+1$.
Thus
$$
 \langle Mu,\,1\rangle=\sum
 \,u_r\,(r+1)^{-2}.
$$
Now take in (\ref{getridlog}) $\,f=p_m$, $\,g=p_n$,
hence $\,u=p_m\,q_{n+1}\,$. Replacing $\,u_{r+1}\,$
by $\,S(m,n,r)$,
we arrive at (\ref{ilog2sum}).
\WhiteBox

The scalar product $\,\langle p_m,\,Np_n\rangle\,$
can be found without difficulty, if $\,m\geq n$.

\begin{Lemma}
\label{evalNpp}$\;$
{\rm (i)} $\,$ If $\;m\geq n+2\,$, then
$\;\langle p_m,\,Np_n\rangle=0$.
\\[\smallskipamount]
{\rm (ii)} $\,$ Nonvanishing scalar products with $\,m\geq n$ are
\begin{equation}
\label{diagNpp}
  \langle p_n,\,Np_n\rangle=[(2n+1)(2n+2)]^{-1},
\end{equation}
\begin{equation}
\label{Npp1}
  \langle p_{n+1},\,Np_n\rangle=\frac{\hnum{n+1}}%
  {2\,(2n+1)\,(2n+3)},
\end{equation}
where $\,\hnum{r}\,$ denotes the $\,r$-th harmonic sum
\begin{equation}
\label{hsum}
 \hnum{r}\;=\;1\,+\,{2}^{-1}\,+\dots+\,{r}^{-1}.
\end{equation}
\end{Lemma}

\Proof$\;$
(i) follows from orthogonality, since $\,\deg Np_n =n+1<m$.

\noindent
(ii)$\,$ Apply the operator $N$ to a Maclaurin series
(\ref{genericMac}). If $\,U=\dst\int u$, then
$$
  \langle u\rangle_t^x \,=\,\frac{U(x)-U(t)}{x-t}\,=\,
  \sum_{r\geq 0} \frac{u_r}{r+1}\,\sum_{j=0}^r x^{r-j} \,t^j.
$$

\noindent
Therefore,
\begin{equation}
\label{NuMac}
 Nu(x)=\sum_{r\geq 0} u_r\,\frac{h_{r+1}}{r+1}\,x^{r+1}.
\end{equation}
Considering the expansion of $\,Np_n$,
we obtain (in notation (\ref{pMac}))
$$
\begin{array}{l}
 Np_n(x)\;=\;\dst a_n^{(n)}\,\frac{h_{n+1}}{n+1}\,x^{n+1}\,+\,
 a_{n-1}^{(n)}\,\frac{h_{n}}{n}\,x^{n}
 \\[3ex] \qquad\quad =\;\dst
 \frac{a_{n}^{(n)}}{a_{n+1}^{(n+1)}}\,\frac{h_{n+1}}{n+1}\,p_{n+1}(x)
 \,+\,
 \left(\,
 \frac{a_{n-1}^{(n)}}{a_{n}^{(n)}}\,\frac{h_{n}}{n}\,-\,
 \frac{a_{n}^{(n+1)}}{a_{n+1}^{(n+1)}}\,\frac{h_{n+1}}{n+1}
 \,\right)\,p_n(x)+\dots.
\end{array}
$$
Since
$
 \; \dst n^{-1}\,a_{n-1}^{(n)}/a_{n}^{(n)}=-1/2
 \;
$
is independent of $\,n$, the coefficient at $\,p_n\,$ in $\,Np_n\,$
equals
$\;(-1/2)\,(h_{n}-h_{n+1})= (2n+2)^{-1}\,$.
Recalling the value (\ref{pL2}) of $\,\langle p_n,p_n\rangle$,
we obtain (\ref{diagNpp}).
The coefficient at $\,p_{n+1}\,$ in $\,Np_n\,$
and consequently the value (\ref{Npp1}) also follow by (\ref{pMac})
and (\ref{pL2}).
\WhiteBox

\subsection{Recurrence for $\;S(m,n,r)$.}

\begin{Lemma}
The following expression
is symmetric in $\,m$ and $n\;$ $\,(m,n\geq 1)$
\begin{equation}
\label{defRmnp}
\Rs(m,n,r)\,=\,(n+1)\;S(m,n,r)\;-\;(n-1)\;S(m-1,n-1,r)
\end{equation}
\end{Lemma}
\Proof$\;$
The generating function $\,\Rs_{mn}(x)\,$ for the quantities
$\,\Rs(m,n,\cdot)\,$ is
$$
\sum_{r=0}^{m+n} \Rs(m,n,r)\,x^{r+1} \,=\, (n+1)\,
p_m(x)\,q_{n+1}(x) \,-\,(n-1)\,p_{m-1}(x)\,q_{n}(x).
$$
By (\ref{q2p}),
$$
\Rs_{mn}\,=\,\frac{n+1}{2n+1}\,p_m\,(p_{n+1}-p_{n-1})\,-
\,\frac{n-1}{2n-1}\,p_{m-1}\,(p_{n}-p_{n-2}).
$$
Apply (\ref{pRec}) to $\,p_{n+1}\,$ and
$\,p_{n}\,$ and find
$$
\Rs_{mn}(x)\,=\,\tilde x\,p_m\,p_n \,-\,p_m\,p_{n-1}
\,-\,\frac{n-1}{n} \,\tilde x\,p_{m-1}\,p_{n-1}
\,+\,\frac{n-1}{n} \,p_{m-1}\,p_{n-2}.
$$
The first term is explicitly symmetric. Apply (\ref{pRec})
to $\,p_m\,$ in the 2nd term and obtain a symmetric formula
for $\;\Rs_{mn}(x)\,-\,\tilde x\,p_m\,p_n \,$:
$$
\left(1+\frac{n-1}{n}+\frac{m-1}{m}\right)\tilde x\,p_{m-1}p_{n-1}
\,+\,\left(\frac{m-1}{m} \,p_{m-2} p_{n-1}+
\frac{n-1}{n} \,p_{m-1} p_{n-2}\right).
$$
\WhiteBox

\noindent{\bf Corollary.}
$\,\DbL(m,n)\,$ satisfy the recurrence relation
\begin{equation}
\label{recLmn}
\begin{array}{l}
\dst (n+1)\,\DbL(m,n)\,-\,(m+1)\,\DbL(n,m) \\[2ex]
\qquad\quad=\;\dst
(n-1)\,\DbL(m-1,n-1)\,-\,(m-1)\,\DbL(n-1,m-1).
\end{array}
\end{equation}

\subsection{Completion of proof}

\begin{Lemma}
\label{Dblsymm}
$\;$
{\rm (i)} $\,$ If $\;|m-n|\geq 2\,$ then
$\;\DbL(m,n)+\DbL(n,m)=0$.
\\[\smallskipamount]
{\rm (ii)} $\,$ In the cases when $\,|m-n|$, we have
\begin{equation}
\label{diagL}
  2\,\DbL(n,n)\,=[n(2n+1)(2n+2)]^{-1}.
\end{equation}
\begin{equation}
\label{Ld1}
  \DbL(n+1,n)\,+\,\DbL(n,n+1)\,=\,-[(2n+1)(2n+2)(2n+3)]^{-1}.
\end{equation}
\end{Lemma}

\Proof$\;$
From proof of Lemma \ref{logred} we see that
$$
 \DbL(m,n)+\DbL(n,m)\,=\,\langle M[p_m\,q_{n+1}+q_{m+1}\,p_n],\,1\rangle
 \,=\,\int\limits_0^1 \frac{q_{m+1}(x)\,q_{n+1}(x)}{x}\,dx.
$$
If $\,n<m-1\,$ then $\,q_{m+1}$ is orthogonal to any polynomial
of degree $n\,$ (see (\ref{q2p})), in particular to $\,x^{-1} q_{n+1}(x)$.
This proves (i).

\noindent
(ii)$\,$
Set
$$
 w_{n+1}(x)=\int_0^x x^{-1}\,q_{n+1}(x)\,dx.
$$
Integrating by parts and remembering that $\,q_{n+1}(0)=q_{n+1}(1)=0\,$
for $\,n>0$, we get
$$
 \DbL(m,n)+\DbL(n,m)\,=\,-\int_0^1 p_m(x)\,w_{n+1}(x)\,dx.
$$
Find the $\,(n+1)$-th and $\,n$-th Legendre coefficients
of $\,w_{n+1}\,$ by considering the respective Maclaurin
coefficients. The calculation is very similar to that in proof of
Lemma \ref{evalNpp}
$$
\begin{array}{l}
 w_{n+1}(x)\;=\;\dst \frac{a_n^{(n)}}{(n+1)^2}\,x^{n+1}\,+\,
 \frac{a_{n-1}^{(n)}}{n^2}\,x^{n}
 \\[3ex] \qquad\quad=\;\dst
 \frac{a_{n}^{(n)}}{a_{n+1}^{(n+1)}}\,\frac{p_{n+1}(x)}{(n+1)^2}
 \,+\,
 \left(\,
 \frac{a_{n-1}^{(n)}}{a_{n}^{(n)}}\,\frac{1}{n^2}\,-\,
 \frac{a_{n}^{(n+1)}}{a_{n+1}^{(n+1)}}\,\frac{1}{(n+1)^2}
 \,\right)\,p_n(x)+\dots.
\end{array}
$$
The coefficient at $\,p_n\,$ equals $\;(-1/2)\,[n^{-1}-(n+1)^{-1}]
=-[2n(n+1)]^{-1}$. The coefficient at $\,p_{n+1}\,$ is
$\;[(2n+1)(2n+2)]^{-1}\,$.
The proof is completed by using (\ref{pL2}),
as in Lemma \ref{evalNpp}
\WhiteBox

\begin{Lemma}
\label{Lemmzero}
$\;$
$\,$ If $\,m\geq 2$, then
\begin{equation}
\label{Lmzero}
  \DbL(m,0)\,=\,(-1)^{m+1}\,[m(m+2)(m^2-1)]^{-1}.
\end{equation}
\end{Lemma}

\Proof$\;$
By (\ref{pRec}),
$$
 2\,p_m(x)\,q_{1}(x)\equiv\, (\tilde x + 1)\,p_m(x) =\,
 \frac{m+1}{2m+1}\,p_{m+1}(x)\,+\,p_m(x)\,+\,\frac{m}{2m+1}\,p_{m-1}(x).
$$
Applying the operator $M$ and integrating, we get
$$
 2\,\DbL(m,0)\,=\,\int\limits_0^1\,
 \left(
 \frac{m+1}{2m+1}\,\frac{q_{m+2}(x)}{x}\,+
 \,\frac{q_{m+1}(x)}{x}\,+\,\frac{m}{2m+1}\,\frac{q_{m}(x)}{x}
 \right)\,dx.
$$
Now use (\ref{qx2p}), and (\ref{Lmzero}) follows .
\WhiteBox

\noindent
Final evaluation of $\,\Hlog(m,n)\,$ requires separate
treatment of the cases $\,m-n\geq2$, $\,m=n$, and $\,m-n=1$.

\medskip
\noindent
\underline{Case I}: $\;m-n\geq2$.
In this case, $\,\DbL(m,n)=-\DbL(n,m)\,$ by Lemma \ref{Dblsymm} (i), and
(\ref{recLmn}) becomes
$$
 (m+n+2)\, \DbL(m,n)\,=\,(m+n-2)\, \DbL(m-1,n-1).
$$
The value $\,\DbL(m,0)\,$ is known from (\ref{Lmzero}).
By induction we find
$$
 \DbL(m,n)=\,-[(m+n)(m+n+2)((m-n)^2)]^{-1}.
$$
Since the second term in the r.h.s.\ of (\ref{ilog2sum}) vanishes
--- see Lemma \ref{evalNpp} (i) --- we obtain the answer (\ref{logleghs}).

\medskip
\noindent
\underline{Case II}: $\;m=n$.
Using (\ref{diagL}) and
(\ref{diagNpp}), we find the value of (\ref{ilog2sum}):
$$
  -\,\frac{1}{2n\,(2n+1)\,(2n+2)}\,-\,
  \frac{1}{(2n+1)\,(2n+2)}\;=\,\frac{1}{2n\,(2n+2)(-1)},
$$
which agrees with (\ref{logleghs}).

\medskip
\noindent
\underline{Case III}: $\;m=n+1$.
Denote
$\;\DbL(n+1,n)=\lambda(n)$.
The recurrence (\ref{recLmn}) with the nonhomogeneous symmetry relation
(\ref{Ld1}) leads to the nonhomogeneous linear difference equation
$$
(2n+3)\lambda_n - (2n-1)\lambda_{n-1}=\frac{2n+5}{(2n-1)\,
(2n+1)_{3}}.
$$
A solution of the homogeneous part is
$\;\lambda^{(0)}_n=[(2n+1)(2n+3)]^{-1}$. Set
$\;\lambda_n=\lambda^{(0)}_n\,\mu_n$. Then $\,\mu_n\,$ satisfies
the equation
$$
\mu_n-\mu_{n-1}=\frac{2n+5}{(2n-1)(2n+2)_2}\,
=\,\frac{1}{2}\,\left(\frac{1}{2n-1}\,+\,\frac{1}{2n+3}\,-\,
\frac{1}{n+1}
\right).
$$
By telescoping,
\begin{equation}
\label{musum}
 \mu_n\,=\mu_0-\,\frac16\,+\,\hodd{2n+1}\,-\,\frac{1}{2}\,\hnum{n+1}\,
 +\frac{1}{2}\left((2n+3)^{-1}\,-\, (2n+1)^{-1}\right).
\end{equation}
Explicit calculation shows that $\mu_0=3K(1,0)=-1/12$.
Substituting (\ref{musum}) together with
(\ref{Npp1}) to (\ref{ilog2sum}), we obtain (\ref{logleg1s}).

\noindent
The proof is complete.
\BlackBox

\smallskip\noindent
{\bf Excercise.}
Let $\,m$, $n$, $l\,$ be nonnegative integers.
Evaluate the 
sum: 
$$
T(m,n,l)=
\sum_{p=0}^{m+n} \frac{S(m,n,p)}{p+l}.
$$
%
\appendix
\section{Legendre polynomials $\,P_n(x)$}
\ubsubsection{Definition.}
$\,P_n(x)\,$ is the polynomial of degree $n$
such that: $\,$ (i) $\,P_n(1)=1$, and $\,$ (ii)
for any polynomial $\,f(x)\,$ of degree $<n$
$$
\langle f,\,P_n\rangle\,\equiv\,\int\limits_{-1}^{1} f(x)\, P_n(x)\,dx\,=\,0.
$$
\ubsubsection{Rodrigues formula}
\begin{equation}
\label{PRod}
 P_n(x)\,=\,\frac{1}{2^n n!}\;\frac{d^n}{dx^n}\, (x^2-1)^n.
\end{equation}
%
\ubsubsection{Symmetry}
\begin{equation}
\label{PSym}
 P_n(-x)=(-1)^n P_n(x).
\end{equation}
\ubsubsection{$\mathbf{L_2}$ norm}
\begin{equation}
\label{PL2}
\langle P_n,\,P_n\rangle\,=\,
2\,(2n+1)^{-1}.
\end{equation}
\ubsubsection{Recurrence relations}
\begin{equation}
\label{PRec}
 (n+1)\,P_{n+1}(x)\,=\,(2n+1)\,x\,P_n(x)\, -\, n\,P_{n-1}(x),
\end{equation}
\begin{equation}
\label{PRecd}
P_{n+1}'(x)\,=\,(2n+1)\,P_n(x)\,+\,P_{n-1}'(x).
\end{equation}

\noindent
\ubsubsection{Shifted Legendre polynomials}$\;$\cite[22.2.11]{AS}
\begin{equation}
\label{pdef}
p_n(x)=P_n(\tilde x),\qquad \tilde x=2x-1,
\end{equation}
form an orthogonal basis on $[0,1]$. When working with $p_n$'s,
we denote by $\langle\cdot,\cdot\rangle$ the integral scalar product
on $[0,1]$.
For reference, here are the shifted versions of the above formulae:
\begin{equation}
\label{pRod}
p_n(x)=(n!)^{-1}\;(d/dx)^n\,\left(\,x^n(x-1)^n\right)\,,
\end{equation}
\begin{equation}
\label{pL2}
\langle p_n,\,p_n\rangle\,=\,
(2n+1)^{-1}\,,
\end{equation}
\begin{equation}
\label{pRec}
 (n+1)\,p_{n+1}(x)\,=\,(2n+1)\,\tilde x\,p_n(x)\, -\, n\,p_{n-1}(x).
\end{equation}

\noindent
\ubsubsection{A bilinear identity for derivatives.}$\;$
If $\,r=|m-n|+2s$, $\,s\geq 0$, then
\begin{equation}
\label{pMacid}
\sum_{j+k=r+1}\, (-1)^j\,
p_m^{(j)}(0)\, p_n^{(j)}(0)\, =\,0.
\end{equation}
\Proof$\;$
The same combination of derivatives taken at $1$ has opposite sign.
So, by integration by parts, the l.h.s.\ equals to
$$
(\pm 1/2)\,(\,
\langle p_m,\, p_n^{(r+1)}\rangle \,-\,
\langle p_m^{(r+1)},\,p_n\rangle\,),
$$
which is $0$, because, e.g.,
$\;\deg (p_m^{(r+1)})=\max(0,\,n-1-2s)\, < n$.
\WhiteBox

\noindent
\ubsubsection{Maclaurin expansion} 
\begin{equation}
\label{pMac}
 p_n(x)=\sum_{k=0}^n  a_{k}^{(n)} x^k,
 \qquad
 a_{k}^{(n)}= \frac{(-1)^{n-k}\,(n+k)!}{(n-k)!\,k!^2}.
\end{equation}

\noindent
\ubsubsection{Integrated shifted Legendre polynomials $\,q_{n+1}(x)$}
\begin{equation}
\label{ipMac}
 q_{n+1}(x)\,=\,\int_0^x p_n(t)\,dt\,=\,
 \sum_{k=0}^n  \frac{a_{k}^{(n)}}{k+1}\, x^{k+1}.
\end{equation}
These polynomials are used in Sect.~5. The notation is local to this work.

\noindent
\ubsubsection{Decomposition} $\;$ (Cf. (\ref{PRecd}).)
$\,$ For $n>0$,
\begin{equation}
\label{q2p}
 2\,(2n+1)\,q_{n+1}\,=\,p_{n+1}\,-\,p_{n-1}.
\end{equation}
\ubsubsection{Integral $\;\langle q_{n+1},x^{-1} \rangle$.}
$\;$ If $n>0$, then
\begin{equation}
\label{qx2p}
 \int\limits_0^1
 \frac{q_{n+1}(x)}{x}\,dx \,\equiv \,
 -\int\limits_0^1
 p_{n}(x)\,\ln{x}\,dx  \,=\,
 \frac{(-1)^{n}}{n(n+1)}.
\end{equation}
\Proof$\;$
 Set $\,w(x)=x(x-1)$. By the Rodrigues formula
 $\,q_{n+1}=(n!)^{-1} \partial^{n-1} w^n$.
 Integrating by parts $n$ times, we get
 $$
 \int\limits_{\eps}^1\frac{q_{n+1}}{x}\,dx \,=\,
 \sum_{j=0}^{n-2}\,\frac{j!}{n!}\,\left[
 \frac{\partial^{n-2-j}\, w^n}{x^{j+1}}
 \right]_{\eps}^1\;
  +\;\frac{(n-1)!}{n!}\,\int\limits_{\eps}^1\frac{w^n}{x^n}\,dx.
 $$
 When $\eps\to0^{+}$, the only survivor is
 $\; n^{-1}\int\limits_{0}^1 (x-1)^n\,dx=(-1)^{n}[n(n+1)]^{-1}
 $. \WhiteBox

\section{Hypergeometric series and identities}
%
%
\ubsubsection{Rising factorial (Pohhammer's symbol).}$\;$
Let $n\in \Zz$, $a\in \Cc$. By definition,
$\; (a)_n\,=\,{\Gamma(a+n)}/{\Gamma(a)}$.
In particular, $\,(1)_n=n!$.
Alternative definition:
\begin{equation}
\label{Pohh}
\begin{array}{ll}
{\rm(i)}\quad\; &
 (a)_0\,=\,1\,;
\\[1.2ex]
{\rm(ii)}\quad\; &
 (a)_n\,=\,a(a+1)\cdot\dots\cdot(a+n-1)
 \quad\;\mbox{\rm if}\;\;n>0;
\\[1.2ex]
{\rm(iii)}\quad &
 (a)_{n}\,=\,[(a+n)_{-n}]^{-1}
 \quad\;\mbox{\rm if}\;\;n<0.
\end{array}
\end{equation}
Three useful formulae:
\begin{equation}
\label{Pohsum}
 (a)_{m+n}=(a)_m\,(a+m)_n,
\end{equation}
\begin{equation}
\label{Pohinv}
(a)_n= (-1)^n\,(-a-n+1)_{n},
\end{equation}
\begin{equation}
\label{Pohdup}
(a)_n\,(a+1/2)_n= 4^{-n}\,(2a)_{2n},
\end{equation}
\ubsubsection{Hypergeometric series.} 
$\;$
The series depending on $p+q+1$ parameters%
\begin{equation}
\label{hgdef}
_p F_q\left(
  {{a_1,\dots,a_p} \atop {b_1,\dots,b_q}}
  \,;\, z
  \right) \;=\;
  \sum_{n=0}^{\infty}
  \frac{(a_1)_n\dots (a_p)_n}{n!\,(b_1)_n\dots (b_q)_n}\,
 z^n.
\end{equation}
is called the (generalized) hypergeometric series.
The ordinary (Gauss') series
is $_2F_1$.
If one of the $\,a_j\,$'s in (\ref{hgdef})
is a negative integer $-r$, then the series
turns to a polynomial of degree $r$ in $z$ and is called {\bf terminating}.
In the frequent case $\mathbf{\,z=1}$, the argument $z$
of $\,_pF_q\,$ is omitted.

%

\smallskip
\noindent
\ubsubsection{Euler's integral for $\,_2F_1$ \quad
{\normalfont \cite[(15.3.1)]{AS}, \cite[2.1.3 (10)]{Bat},
\cite[1.5]{Bai}}.
}
\\
If $\;\Re c>\Re b>0\;$ and $\,z\notin[1,\infty)$, then
\begin{equation}
\label{f21int}
 \HGa{a}{b}{c}{z}=
 \int\limits_0^1 t^{b-1}\,(1-t)^{c-b-1}\,(1-tz)^{-a}\,dt .
\end{equation}
\ubsubsection{Gauss' formula for $\,_2F_1(\dots;\,1)$ \quad
{\normalfont \cite[(15.1.20)]{AS}, \cite[2.1.3 (14)]{Bat},
\cite[1.3]{Bai}}
}
\begin{equation}
\label{hgGauss}
 \HGu{a}{b}{c}\,=\,\frac{\Gamma(c)\,\Gamma(c-a-b)}%
 {\Gamma(c-a)\,\Gamma(c-b)} \qquad
(\Re (c-a-b)>0).
\end{equation}
\ubsubsection{Thomae's transformation formula for $\,_3F_2$
$\quad${\normalfont \cite[(2)]{W}, \cite[3.2]{Bai},
\cite[(3.1.2)]{GR}}
}
If $s=e+f-a-b-c$, then
\begin{equation}
\label{thomae32}
 \HGb{a}{b}{c}{e}{f}=
 \frac{\Gamma(e)\,\Gamma(f)\,\Gamma(s)}{\Gamma(a)\,\Gamma(s+b)
 \,\Gamma(s+c)}\,
 \HGb{e-a}{f-a}{s}{s+b}{s+c}.
\end{equation}

\noindent
\ubsubsection{Whipple's summation formula for $\,_3F_2$
$\;${\normalfont \cite[(34)]{W}, \cite[3.4]{Bai},
\cite[4.4 (7)]{Bat}}
}.
Terminating case: for a nonnegative integer $n$,
\begin{equation}
\label{whipterm}
\HGb{-n}{n+1}{(b+c-1)/2}{b}{c}
\;=\;
4^n\,
\frac{(\frac{b-n}{2})_n\, (\frac{c-n}{2})_n}
{(b)_n (c)_n}.
\end{equation}
General case: let $\,a\in\Cc$,
$\,a'=1-a$, and $\,\Re(b+c)>1$. Then
\begin{equation}
\label{whip}
\HGb{a}{a'}{(b+c-1)/2}{b}{c}\;=\;
\frac{4\pi\,\cdot 2^{-(b+c)}\, \Gamma(b)\,\Gamma(c)}
{
\Gamma(\frac{b+a}{2})\, \Gamma(\frac{b+a'}{2})\,
\Gamma(\frac{c+a}{2})\, \Gamma(\frac{c+a'}{2}).
}
\end{equation}
\noindent
\ubsubsection{Dougall's formula for $\,_5F_4$
$\quad${\normalfont \cite[4.4]{Bai}, \cite[4.5 (6)]{Bat}}
}.
Suppose that
$$
1+a\,=\,b+b'\,=\,c+c'\,=\,d+d'\,=\,e+e',
$$
and $\,b-b'=1$. Then
\begin{equation}
\label{dougall54}
 \HGc{a}{b}{c}{d}{e}{b'}{c'}{d'}{e'}
 \;=\;
 \frac{\Gamma(c')\,\Gamma(d')\,\Gamma(e')\,\Gamma(1+a-c-d-e)}%
 {\Gamma(1+a)\,\Gamma(c'-d)\,\Gamma(d'-e)\,\Gamma(e'-c)}.
\end{equation}
\ubsubsection{Particular $\,_3F_2\,$ case of Karlsson's reduction formula}
$\;${\normalfont \cite{K}}
\begin{equation}
\label{Karlsson}
 \begin{array}{l}
 \dst\HGb{-b-n}{b+n+1}{d}{b}{c}
\\[2ex]
\dst
 \qquad\quad
 =\;\sum_{k=0}^{n+1}\Binomial{n+1}{k}
 \frac{(-b-n)_k\, (d)_k}{(b)_k\, (c)_k}
 \HGu{-b-n+k}{d+k}{c+k}.
 \end{array}
\end{equation}



%
\end{document}